\documentclass[11pt]{amsart}

\newtheorem{theorem}{Theorem}[section]
\newtheorem{proposition}[theorem]{Proposition}
\newtheorem{lemma}[theorem]{Lemma}

\newcommand{\abs}[1]{\lvert#1\rvert}   
\newcommand{\set}[1]{\{ #1 \}}         

\newcommand{\Ab}{\operatorname{Ab}}
\newcommand{\Aut}{\operatorname{Aut}}
\newcommand{\cent}{\operatorname{cent}}
\renewcommand{\deg}{\operatorname{deg}}
\newcommand{\Diff}{\operatorname{Diff}}
\newcommand{\Dih}{\operatorname{Dih}}
\newcommand{\GL}{\operatorname{GL}}
\newcommand{\Isom}{\operatorname{Isom}}
\newcommand{\isom}{\operatorname{isom}}
\newcommand{\Norm}{\operatorname{Norm}}
\newcommand{\norm}{\operatorname{norm}}
\newcommand{\Or}{\operatorname{O}}
\newcommand{\Out}{\operatorname{Out}}
\newcommand{\SO}{\operatorname{SO}}

\newcommand{\C}{\operatorname{{\mathbb C}}}
\renewcommand{\P}{\operatorname{{\mathbb P}}}
\newcommand{\R}{\operatorname{{\mathbb R}}}
\newcommand{\Z}{\operatorname{{\mathbb Z}}}

\renewcommand{\H}{\operatorname{{\mathcal H}}}
\newcommand{\Is}{\operatorname{{\mathcal I}}}
\newcommand{\orb}{\operatorname{{\mathcal O}}}

\newcommand{\Ostar}{\Or(2)^*}
\newcommand{\ttimes}{\operatorname{\,\widetilde{\times}\,}}

\begin{document}

\title{Isometries of elliptic $3$-manifolds}

\author{Darryl McCullough}
\address{Department of Mathematics\\
University of Oklahoma\\
Norman, Oklahoma 73019\\
USA}
\email{dmccullough@math.ou.edu}
\urladdr{www.math.ou.edu/$_{\widetilde{\phantom{n}}}$dmccullo/}

\date{October 1, 2000}


\begin{abstract}
The closed $3$-manifolds of constant positive curvature were classified
long ago by Seifert and Threlfall. Using well-known information about the
orthogonal group $\Or(4)$, we calculate their full isometry groups
$\Isom(M)$, determine which elliptic $3$-manifolds admit Seifert fiberings
that are invariant under all isometries, and verify that the inclusion of
$\Isom(M)$ to $\Diff(M)$ is a bijection on path components.
\end{abstract}

\maketitle
\markboth{DARRYL McCULLOUGH}{ISOMETRIES OF $3$-MANIFOLDS}

\section*{Introduction}

A closed $3$-manifold with a Riemannian metric of constant curvature~$1$ is
called an elliptic $3$-manifold. These were classified long ago by Seifert
and Threlfall. Using well-known information about the group $\Or(4)$ of
isometries of $S^3$, it is a routine task, albeit a laborious one, to
calculate the isometry groups $\Isom(M)$ for all such $M$.  This may be
``folk'' knowledge, indeed it is a logical continuation of the beautiful
expositions in~\cite{Scott} and~\cite{Sakuma}. But we are not aware of any
explicit reference for this calculation, and we believe that details should
appear in the literature.  These comprise most of this paper, and the
results are summarized in tables~\ref{tab:isometries} and~\ref{tab:lens
spaces} below.

A calculation of $\Isom(M)$ for lens spaces (the most difficult case)
appears in a paper by J. Kalliongis and A. Miller \cite{KM}. The authors
apply their calculation to describe all finite group actions by isometries
on lens spaces.  The methodology of that paper is rather different from
ours. It exploits the genus-$1$ Heegaard splittings of lens spaces to
provide a concrete description of the isometries. In contrast, we use a
very algebraic approach, which applies uniformly to all elliptic
$3$-manifolds, but does not furnish (at least, not without further
examination) such a meaningful description. Also, considerable information
about the isometries of lens spaces was obtained in~\cite{HR}.

After obtaining our explicit calculation of $\Isom(M)$, we use it to check
that all except five of the elliptic $3$-manifolds have Seifert fiberings
preserved by $\isom(M)$, the connected component of the identity, and that
all except six have Seifert fiberings preserved by the group $\Isom_+(M)$
of orientation-preserving isometries. Making use of the calculations of the
mapping class group $\pi_0(\Diff(M))$ by various authors, we verify that
the inclusion from $\Isom(M)$ to the diffeomorphism group $\Diff(M)$
induces a bijection on path components. This is the ``$\pi_0$'' part of the
Generalized Smale Conjecture, which asserts that the inclusion from
$\Isom(M)$ to $\Diff(M)$ is a homotopy equivalence for any elliptic
$3$-manifold. In fact, we show that there is always a subgroup of
$\Isom(M)$ that is carried isomorphically to $\pi_0(\Diff(M))$.  The
Generalized Smale Conjecture is known for $S^3$ \cite{H1}, for all
quaternionic and prism manifolds~\cite{I1,I2,M-R}, and for lens spaces
$L(4k,2k-1)$ with $k\geq 2$~\cite{M-R}.

We appreciate the receipt of a copy of \cite{KM} from one of its
authors. Their correct calculation of $\Isom(M)$ for lens spaces revealed
an error in an earlier version of our work.

\section[Isometries]{Isometries}
\label{sec:isometries intro}

Let $M$ be a closed $3$-manifold with a Riemannian metric of constant
positive curvature. By rescaling, we may assume that $M$ is elliptic. There
is a developing map from the universal cover of $M$ to the standard
$3$-sphere $S^3$, through which $\pi_1(M)$ corresponds to a group $G$ of
isometries acting freely on $S^3$. The full group of isometries of $S^3$ is
the orthogonal group $\Or(4)$, but the Lefschetz fixed-point formula shows
that every orientation-reversing homeomorphism of $S^3$ has a fixed point,
so $G$ must be a subgroup of~$\SO(4)$.  The finite subgroups of $\SO(4)$
that act freely on $S^3$ were determined by Hopf and
Seifert-Threlfall~\cite{S-T1,S-T2}, and reformulated using quaternions by
Hattori~\cite{H}. Among the references for this material are \cite{Wolf}
(pp.~226-227), \cite{Orlik} (pp.~103-113), \cite{Scott} (pp.~449-457),
\cite{P}, and~\cite{Sakuma}. The latter calculates the isometry groups of
certain orbifold quotients of~$S^3$.

When $G$ is not cyclic, any two subgroups of $\SO(4)$ that are isomorphic
to $G$ and act freely are conjugate by an element of $\Or(4)$, so there is
a unique isometry class of elliptic $3$-manifold for each such $G$. For the
cyclic case, the conjugacy classes correspond to the lens spaces $L(m,q)$
with~$q\leq m/2$. So the elliptic structure of $S^3/G$ is uniquely
determined by its homeomorphism type, and in our calculation of $\Isom(M)$
we need only examine one manifold from each homeomorphism type.

In the remainder of this section, we will give a complete calculation of
$\Isom(M)$ for elliptic $3$-manifolds. We first give some preliminary
material on the quaternionic structure of $S^3$, which facilitates explicit
calculations in $\Or(4)$, and on some subgroups of $\SO(4)$ which will
arise.  In section~\ref{sec:calculation of isometries}, we lay out the
general method for calculating $\Isom(M)$.  Section~\ref{sec:isometries}
gives the rather straightforward calculations of $\Isom(M)$ for $M$ other
than lens spaces.  These are summarized in table~\ref{tab:isometries}. The
case of lens spaces is more complicated, and comprises
section~\ref{sec:lens spaces}. Its results are summarized in
table~\ref{tab:lens spaces}. In section~\ref{sec:isometric realization}, we
observe that the group of path components of $\Isom(M)$ can be realized by
a finite subgroup of~$\Isom(M)$.

\subsection{Quaternions and orthogonal transformations}
The $3$-sphere has a rich algebraic structure as the group of unit
quaternions. We begin by reviewing some facts about this structure that
will be used in our later work.  The references~\cite{Scott}
and~\cite{Sakuma} should suffice for all unsupported assertions given in
this section.

Fix coordinates on $S^3$ as $\set{(z_0,z_1)\;|\;z_i\in\C, z_0\overline{z_0}
+ z_1\overline{z_1}=1}$. Its group structure as the unit quaternions can
then be given by writing points in the form $z=z_0+z_1j$, where $j^2=-1$
and $jz_i=\overline{z_i}j$.  The real part $\Re(z)$ is $\Re(z_0)$, and the
imaginary part $\Im(z)$ is $\Im(z_0)+z_1j$.  The inverse of $z$ is
$\Re(z)-\Im(z)=\overline{z_0}-z_1j$.

The unique element of order $2$ in $S^3$ is $-1$, and it generates the
center of $S^3$. The pure imaginary unit quaternions form the $2$-sphere
$P\subset S^3$ of vectors orthogonal to~$1$.

Two elements of $S^3$ are conjugate if and only if they have the same real
part. In particular, $P$ is a single conjugacy class. It is exactly the set
of elements of order~$4$.

Each $z\neq \pm1$ is contained in a unique connected $1$-dimensional
subgroup consisting of all elements of the form $\cos(\theta)+
\sin(\theta)\Im(z)/\|\Im(z)\|$, which is exactly the centralizer of $z$.
Any two such subgroups are conjugate.

By $S^1$ we will denote the subgroup of points in $S^3$ with $j=0$, that
is, all $z_0\in S^1\subset \C$.  Let $\xi_k=\exp(2\pi i/k)$, which
generates a cyclic subgroup $C_k\subset S^1$. The elements $S^1\cup S^1j$
form a subgroup $\Ostar\subset S^3$, which is exactly the normalizer of
$S^1$ and of the $C_k$ with $k>2$.  When $k$ is odd, the quotient
$\Ostar/C_k$ is isomorphic to $\Ostar$.  When $k$ is even, the quotient is
isomorphic to the orthogonal group $\Or(2)$, which we may describe as
$S^1\cup S^1j_0$ where $j_0^2=1$ and $j_0z_0=\overline{z_0}j_0$ for $z_0\in
S^1$.

The finite subgroups of $S^3$ are well-known: the binary icosahedral,
octahedral, and tetrahedral groups $I_{120}^*$, $O_{48}^*$, and $T_{24}^*$,
the binary dihedral groups $D_{4m}^*$, $m>2$, the finite quaternion group
$Q_8$, and the cyclic groups. One may treat $Q_8$ as a binary dihedral
group $D_8^*$, but we prefer to regard is as a separate family of subgroups
since its normalizer is not just $D_{16}^*$. Two finite subgroups are
conjugate if and only if they are isomorphic. Representatives of each
conjugacy class are given in table~\ref{tab:finite subgroups}, along with
presentations, explicit generators, and information about their normalizers
$\textrm{N}(G)$ and the quotients $\textrm{N}(G)/G$. In our notation, the
subscript attached to the name of a finite group tells its order, except
for $S_3$, the symmetric group on three letters.

\begin{table}
\begin{small}
\renewcommand{\arraystretch}{1.5}
\setlength{\fboxsep}{0pt}
\setlength{\tabcolsep}{4.5pt}
\fbox{%
\begin{tabular}{l|l|l|c|c}
$G$&Presentation&Generators&$\textrm{N}(G)$&$\textrm{N}(G)/G$\\
\hline
\hline
\newlength{\ymwidth}%
$I_{120}^*$&\settowidth{\ymwidth}{$y^m$}$\langle\; x,y\;\vert\;x^2=\makebox[\ymwidth]{$y^3$}=(xy)^5\;\rangle$&%
\newlength{\minipagewidth}%
\settowidth{\minipagewidth}{$y=\frac{1}{2}+\frac{1}{2}\sqrt{1-2\cos(\frac{2\pi}{5})}\,i$}%
\begin{minipage}{\minipagewidth}%
\noindent$x=j$,\rule[2.5 ex]{0mm}{0mm}\par%
\noindent$y=\frac{1}{2}+\frac{1}{2}\sqrt{1-2\cos(\frac{2\pi}{5})}\,i$\par
\noindent$\phantom{y=}+\cos(\frac{\pi}{5})j$\rule[-1.8 ex]{0mm}{0mm}%
\end{minipage}&
$I_{120}^*$&$\set{1}$\\
\hline
$O_{48}^*$&\settowidth{\ymwidth}{$y^m$}$\langle\; x,y\;\vert\;x^2=\makebox[\ymwidth]{$y^3$}=(xy)^4\;\rangle$&%
\settowidth{\minipagewidth}{$y=\frac{1}{\sqrt{2}}\xi_8+\frac{1}{\sqrt{2}}\xi_8j$}%
\begin{minipage}{\minipagewidth}%
\noindent$x=\frac{1}{\sqrt{2}}i+\frac{1}{\sqrt{2}}j$,%
\rule[2.5 ex]{0mm}{0mm}\rule[-1.8 ex]{0mm}{0mm}\par
\noindent$y=\frac{1}{\sqrt{2}}\xi_8+\frac{1}{\sqrt{2}}\xi_8j$\rule[-1.8 ex]{0mm}{0mm}
\end{minipage}&
$O_{48}^*$&$\set{1}$\\
\hline
$T_{24}^*$&\settowidth{\ymwidth}{$y^m$}$\langle\; x,y\;\vert\;x^2=\makebox[\ymwidth]{$y^3$}=(xy)^3\;\rangle$&%
$x=j$,\,\;$y=\frac{1}{\sqrt{2}}\xi_8+\frac{1}{\sqrt{2}}\xi_8j$\rule[-1.8 ex]{0mm}{0mm}&%
$O_{48}^*$&$C_2$\\
\hline
$D_{4m}^*$&$\langle\; x,y\;\vert\;x^2=y^m=(xy)^2\;\rangle$&%
$x=j$,\,\;$y=\xi_{2m}$&%
$D_{8m}^*$&$C_2$\\
\hline
$Q_8$&\settowidth{\ymwidth}{$y^m$}$\langle\; x,y\;\vert\;x^2=\makebox[\ymwidth]{$y^2$}=(xy)^2\;\rangle$&%
$x=j$,\;\,$y=i$&$O_{48}^*$&$S_3$\\
\hline
$C_{2k-1}$&$\langle\; t\;\vert\;t^{2k-1}=1\;\rangle$&$t=\xi_{2k-1}$&$\Ostar$&$\Ostar$\\
\hline
$C_{2k}$&$\langle\; t\;\vert\;t^{2k}=1\;\rangle$&$t=\xi_{2k}$&$\Ostar$&$\Or(2)$\\
\hline
$C_2$&$\langle\; t\;\vert\;t^2=1\;\rangle$&$t=-1$&$S^3$&$\SO(3)$
\end{tabular}}
\end{small}
\bigskip
\caption{Finite subgroups of $S^3$ ($m>2$, $k>1$)}
\label{tab:finite subgroups}
\end{table}

The usual inner product on $S^3$ is given by $z\cdot
w=\Re(zw^{-1})$. Consequently, left multiplication and right multiplication
by elements of $S^3$ are orthogonal transformations of $S^3$, and there is
a homomorphism $F\colon S^3\times S^3\to \SO(4)$ defined by
$F(q_1,q_2)(q)=q_1 q q_2^{-1}$. It is surjective and has kernel
$\set{(1,1),(-1,-1)}$.  The center of $\SO(4)$ has order~2, and is
generated by $F(1,-1)$, the antipodal map of~$S^3$.

We use the following notational convention. If $H_1$ and $H_2$ are groups,
each containing $-1$ as a central involution, then the quotient $(H_1\times
H_2)/\langle (-1,-1)\rangle$ will be denoted by $H_1\ttimes H_2$. In
particular, $\SO(4)$ itself is $S^3\ttimes S^3$, and contains the subgroups
$\Ostar\ttimes S^3$, $S^1\ttimes S^3$, $\Ostar\ttimes\Ostar$, and
$S^1\ttimes S^1$. The latter is isomorphic to $S^1\times S^1$, but as we
will see, it is sometimes useful to distinguish between them.

There are $2$-fold covering homomorphisms
\begin{equation*}
\Ostar\times \Ostar\to
\Ostar\ttimes\Ostar\to
\Or(2)\times \Or(2)\to
\Or(2)\ttimes\Or(2)\ .
\end{equation*}
Each of these groups is diffeomorphic to four disjoint copies of the
torus. The first three are pairwise nonisomorphic. Indeed, they are
distinguished by their conjugacy classes of involutions, having three,
five, and eight respectively. The fourth group, $\Or(2)\ttimes\Or(2)$, is
isomorphic to $\Ostar\ttimes\Ostar$. An isomorphism from
$\Ostar\ttimes\Ostar$ to $\Or(2)\ttimes\Or(2)$ is defined on generators by
sending $(x_0,x_1)$ to $(x_0,x_1)$ for $x_i\in S^1$, $(j,1)$ to $(j_0,i)$,
and $(1,j)$ to $(i,j_0)$.

Similarly, there are $2$-fold covering homomorphisms of the pairwise
nonisomorphic groups
\begin{equation*}
S^1\times S^3\to
S^1\ttimes S^3\to
S^1\times \SO(3)\ ,
\end{equation*}
although the first two are diffeomorphic.

By $\Dih(S^1\ttimes S^1)$ and $\Dih(S^1\times S^1)$ we denote the
semidirect products $(S^1\ttimes S^1)\circ C_2$ and $(S^1\times S^1)\circ
C_2$, where $C_2$ acts by complex conjugation in both factors. In fact,
these groups are isomorphic, by sending $((z_0,w_0),x)$ to
$((z_0w_0,w_0^2),x)$.  Notice that $S^1\ttimes S^1$ and the element
$F(j,j)$ generate a subgroup $\Dih(S^1\ttimes S^1)$ of index~$2$ in
$\Ostar\ttimes\Ostar$. The element $F(j,1)$ represents its nontrivial
coset.

\subsection{Calculation of $\Isom(M)$}
\label{sec:calculation of isometries}

Since $\Isom(M)$ is a topological group, we may analyze it by determining
$\isom(M)$, the connected component of the identity map, calculating
$\Is(M)$, the group of path components $\pi_0(\Isom(M))$, and understanding
the extension
\begin{equation*}
1\longrightarrow \isom(M)\longrightarrow \Isom(M)\longrightarrow
\Is(M)\longrightarrow 1\ .\tag{$*$}
\end{equation*}
We will see that for elliptic $3$-manifolds, this sequence always splits,
so that $\Isom(M)$ is a semidirect product $\isom(M)\circ \Is(M)$. We use
the notations $\Isom_+(M)$ and $\Is_+(M)$ to indicate the subgroups
consisting of the orientation-preserving elements.

As above, we regard $M$ as the quotient $S^3/G$ where $G$ is a finite
subgroup of $\SO(4)$ acting freely on $S^3$.  Since $\SO(4)=\Isom_+(S^3)$,
$\Isom_+(M)$ is the quotient $\Norm(G)/G$, where $\Norm(G)$ is the
normalizer of $G$ in $\SO(4)$.  Taking connected components of the
identity, $\isom(M)$ is $\norm(G)/(G\cap\norm(G))$.

Assuming that the group $G$ is clear from the context, we denote the
isometry that an element $F(q_1,q_2)$ of $\Norm(G)$ induces on $S^3/G$
by~$f(q_1,q_2)$. By elementary covering space theory, the induced outer
automorphism of $f(q_1,q_2)$ on $\pi_1(M)=G$ equals the outer automorphism
that sends $g$ to $F(q_1,q_2)\,g\,F(q_1,q_2)^{-1}$.

To compute $\Isom_+(M)$, put $G^*= F^{-1}(G)$, and let $G_L$ and $G_R$ be
the projections of $G^*$ into the left and right factors of $S^3\times
S^3$. Since the kernel of $F$ is central, we have $\Norm(G)/G\cong
\Norm(G^*)/G^*$ and $\norm(G)/(G\cap\norm(G))\cong \norm(G^*)/(G^*\cap
\norm(G^*))$.  Since $G^*$ is discrete, $\norm(G^*)$ equals $\cent(G^*)$,
the connected component of the identity in the centralizer of $G^*$. But
$\cent(G^*)=\cent(G_L)\times \cent(G_R)=\norm(G_L)\times \norm(G_R)$, so
$\norm(G^*)=\norm(G_L)\times \norm(G_L)$. We will see, however, that
$\Norm(G^*)$ is sometimes smaller than $\Norm(G_L)\times \Norm(G_R)$. It
will, of course, consist of path components of $\Norm(G_L)\times
\Norm(G_R)$, so to determine which path components lie in $\Norm(G^*)$, we
need only check one representative from each.

\subsection{The case when $M$ is not a lens space}
\label{sec:isometries}

When we speak of lens spaces, we always include the cases of $S^3=L(1,0)$
and $\R\P^3=L(2,1)$.  For the remainder of this section, we assume that $M$
is not a lens space.  By~\cite{N-R} (pp.~188-190, 194), if $M$ is an
elliptic $3$-manifold which is not a lens space, then $M$ does not admit an
orientation-reversing diffeomorphism (or even an orientation-reversing
self-homotopy equivalence). So $\Isom(M)=\Isom_+(M)$.

From \cite{Scott} (see also~\cite{Sakuma}), we have a complete list (up to
conjugacy in $\Or(4)\,$) of the finite subgroups of $\SO(4)$ that act
freely on $S^3$.  In this section we will examine the noncyclic cases.  We
will select our imbeddings of the $G$ into $\SO(4)$ so that the Hopf
fibering of $S^3$ is a maximally symmetric fibering, that is, so that the
subgroup of $\Isom(M)$ that preserves it is as large as possible for any
Seifert fibering on $M$ (see section~\ref{sec:fiberings} for a discussion
of the invariance of Seifert fiberings under $\Isom(M)$). When determining
quotient groups, it is useful to recall from~\cite{Scott} that whenever $G$
is not cyclic, it must contain the antipodal map~$F(1,-1)$.

In Cases I-V, $G^*=G_L\times G_R$, so $\Norm(G^*)=\Norm(G_L)\times
\Norm(G_R)$.
\medskip

\noindent {\sl Case I.} $G=Q_8\times C_n$, $(2,n)=1$
\smallskip

If $n=1$, use $G^*=Q_8\times\set{1}$, and if $n>1$, use $G^*=C_{2n}\times
Q_8$.
\begin{enumerate}
\item
For $n=1$, $\Norm(G^*)\cong O_{48}^*\times S^3$, $\isom(M)\cong \SO(3)$,
and $\Is(M)\cong S_3$ generated by
$f(\frac{1}{\sqrt{2}}i+\frac{1}{\sqrt{2}}j,1)$ and
$f(\frac{1}{\sqrt{2}}\xi_8+\frac{1}{\sqrt{2}}\xi_8j,1)$. These act
trivially on $\isom(M)$, so $\Isom(M)\cong \SO(3)\times S_3$.
\item
For $n>1$, $\Norm(G^*)\cong \Ostar\times O_{48}^*$, $\isom(M)\cong S^1$,
and $\Is(M)\cong C_2\times S_3$ generated by $f(j,1)$,
$f(1,\frac{1}{\sqrt{2}}i+\frac{1}{\sqrt{2}}j)$, and
$f(1,\frac{1}{\sqrt{2}}\xi_8+\frac{1}{\sqrt{2}}\xi_8j)$. Since the second
two generators act trivially on $\isom(M)$, $\Isom(M)\cong \Or(2)\times S_3$.
\end{enumerate}
\smallskip

\noindent {\sl Case II.} $G=D_{4m}^*\times C_n$, $m>2$, $(2m,n)=1$
\smallskip

If $n=1$, use $G^*=D_{4m}^*\times C_2$, and if $n>1$, use $G^*=C_{2n}\times
D_{4m}^*$.
\begin{enumerate}
\item
For $n=1$, $\Norm(G^*)\cong D_{8m}^*\times S^3$, $\isom(M)\cong \SO(3)$,
$\Is(M)\cong C_2$ generated by $f(\xi_{4m},1)$, and $\Isom(M)\cong
\SO(3)\times C_2$.
\item
For $n>1$, $\Norm(G^*)\cong \Ostar\times D_{8m}^*$, $\isom(M)\cong S^1$,
$\Is(M)\cong C_2\times C_2$ generated by $f(j,1)$ and $f(1,\xi_{4m})$, and
$\Isom(M)\cong \Or(2)\times C_2$.
\end{enumerate}
\smallskip

\noindent {\sl Case III.} $G=T_{24}^*\times C_n$, $n\geq1$, $(6,n)=1$
\smallskip

Use $G^*=C_{2n}\times T_{24}^*$.
\begin{enumerate}
\item
If $n=1$, then $\Norm(G^*)\cong S^3\times O_{48}^*$, $\isom(M)\cong
\SO(3)$, $\Is(M)\cong C_2$ generated by
$f(1,\frac{1}{\sqrt{2}}i+\frac{1}{\sqrt{2}}j)$, and $\Isom(M)\cong
\SO(3)\times C_2$.
\item
If $n>1$, then $\Norm(G^*)\cong \Ostar\times O_{48}^*$, $\isom(M)\cong
S^1$, $\Is(M)\cong C_2\times C_2$, generated by $f(j,1)$ and
$f(1,\frac{1}{\sqrt{2}}i+\frac{1}{\sqrt{2}}j)$, and $\Isom(M)\cong
\Or(2)\times C_2$.
\end{enumerate}
\smallskip

\noindent {\sl Case IV.} $G=O_{48}^*\times C_n$, $(6,n)=1$
\smallskip

Use $G^*=C_{2n}\times O_{48}^*$.
\begin{enumerate}
\item
If $n=1$, then $\Norm(G^*)\cong S^3 \times O_{48}^*$,
$\isom(M)=\Isom(M)\cong \SO(3)$, and $\Is(M)=\set{1}$.
\item
If $n>1$, then $\Norm(G^*)\cong \Ostar\times O_{48}^*$, $\isom(M)\cong
S^1$, $\Is(M)\cong C_2$ generated by $f(j,1)$, and $\Isom(M)\cong \Or(2)$.
\end{enumerate}
\smallskip

\noindent {\sl Case V.} $G=I_{120}^*\times C_n$, $(30,n)=1$
\smallskip

Use $G^*=C_{2n}\times I_{120}^*$.
\begin{enumerate}
\item
If $n=1$, then $\Norm(G^*)\cong S^3\times I_{120}^*$,
$\isom(M)=\Isom(M)\cong \SO(3)$, and $\Is(M)\cong \set{1}$.
\item
If $n>1$, then $\Norm(G^*)\cong \Ostar\times O_{48}^*$,
$\isom(M)\cong S^1$, $\Is(M)\cong C_2$ generated by $f(j,1)$, and
$\Isom(M)\cong \Or(2)$.
\end{enumerate}
\medskip

\noindent In the remaining two cases, $G^*$ is a diagonal subgroup of
$G_L\times G_R$. Recall that $H$ is a \emph{diagonal} subgroup of a product
$H_1\times H_2$ when there are surjective homomorphisms $f_i\colon H_i\to
A$ to an abelian group $A$ and $H$ is the kernel of $f_1\times f_2\colon
H_1\times H_2\to A$. This implies that $(H_1\times H_2)/\ker(f_1\times
f_2)\cong H_1/\ker(f_1)\cong H_2/\ker(f_2) \cong A$.
\medskip

\noindent {\sl Case VI.} $G^*$ is a diagonal subgroup of index~$2$ in
$C_{4n}\times D_{4m}^*$, where $m$ is odd and $(m,n)=1$.
\smallskip

We have $\Norm(G_L)\times \Norm(G_R)=\Ostar\times D_{8m}^*$ and
$\isom(M)=S^1$.  Since $m$ is odd, $D_{4m}^*$ has a unique subgroup of
index~$2$, and it follows that $\Norm(G^*)=\Norm(G_L)\times
\Norm(G_R)$. Therefore $\isom(M)=S^1$, $\Is(M)\cong C_2\times C_2$
generated by $f(j,1)$ and $f(1,\xi_{4m})$, and $\Isom(M)\cong \Or(2)\times
C_2$.
\medskip

\noindent {\sl Case VII.}  $G^*$ is a diagonal subgroup of index~3 in
$C_{6n}\times T_{24}^*$, where $n$ is odd and divisible by $3$.\par
\smallskip

We have $\Norm(G_L)\times \Norm(G_R)=\Ostar\times O_{48}^*$ and
$\isom(M)=S^1$.  Since the abelianization of $T_{24}^*$ is $C_3$, there are
exactly two homomorphisms fom $T_{24}^*$ onto $C_3$.  Conjugation by an
element of $O_{48}^*-T_{24}^*$ induces the nontrivial automorphism of the
abelianization.  Similarly, there are two homomorphisms from $C_{6n}$ onto
$C_3$, and conjugation by $j$ induces the nontrivial automorphism on the
quotient. It follows that there are two index~$3$ diagonal subgroups in
$G_L\times G_R$, which are interchanged by $(j,1)$ and by
$(1,\frac{1}{\sqrt{2}}i+\frac{1}{\sqrt{2}}j)$, so $\Norm(G^*)=S^1\times
T_{24}^*\cup (\Ostar-S^1)\times (O_{48}^*-T_{24}^*)$, $\isom(M)=S^1$,
$\Is(M)\cong C_2$ generated
by~$f(j,\frac{1}{\sqrt{2}}i+\frac{1}{\sqrt{2}}j)$, and $\Isom(M)\cong
\Or(2)$.
\medskip

\noindent We note the following containments.
\begin{enumerate}
\item[(a)] In all cases other than $G=T_{24}^*$, $O_{48}^*$, and
$I_{120}^*$, $\norm(G)\subseteq \Ostar\ttimes S^3$.
\item[(b)] In all cases except those in (a) and $G=Q_8$,
$\Norm(G)\subseteq \Ostar\ttimes S^3$.
\end{enumerate}

\noindent Table~\ref{tab:isometries} summarizes our calculations
of~$\Isom(M)$ for the cases other than lens spaces.

\begin{table}
\begin{small}
\renewcommand{\arraystretch}{1.5}
\setlength{\tabcolsep}{2 ex}
\setlength{\fboxsep}{0pt}
\fbox{%
\begin{tabular}{l|l|l|l}
$G$&$M$&$\Isom(M)$&$\Is(M)$\\
\hline
\hline
$Q_8$&quaternionic manifold&$\SO(3)\times S_3$&$S_3$\\ 
\hline 
$Q_8\times C_n$&quaternionic manifold&$\Or(2)\times S_3$&$C_2\times S_3$\\ 
\hline 
$D_{4m}^*$&prism manifold&$\SO(3)\times  C_2$&$C_2$\\ 
\hline 
$D_{4m}^*\times C_n$&prism manifold&$\Or(2)\times C_2$&$C_2\times C_2$\\ 
\hline 
index $2$ diagonal&prism manifold&$\Or(2)\times C_2$&$C_2\times C_2$\\ 
\hline
$T_{24}^*$&tetrahedral manifold&$\SO(3)\times C_2$&$C_2$\\ 
\hline 
$T_{24}^*\times C_n$&tetrahedral manifold&$\Or(2)\times C_2$&$C_2\times C_2$\\ 
\hline 
index $3$ diagonal&tetrahedral manifold&$\Or(2)$&$C_2$\\ 
\hline 
$O_{48}^*$&octahedral manifold&$\SO(3)$&$\set{1}$\\ 
\hline 
$O_{48}^*\times C_n$&octahedral manifold&$\Or(2)$&$C_2$\\ 
\hline 
$I_{120}^*$&icosahedral manifold&$\SO(3)$&$\set{1}$\\ 
\hline 
$I_{120}^*\times C_n$&icosahedral manifold&$\Or(2)$&$C_2$\\
\end{tabular}}
\end{small}
\bigskip
\caption{Isometry groups of $M=S^3/G$\hspace{2 ex}($m>2$, $n>1$)}
\label{tab:isometries}
\end{table}

\subsection{Calculation of $\Isom(M)$ for lens spaces}
\label{sec:lens spaces}

For $(m,q)=1$, the element $F(\xi_{2m}^{q+1},\xi_{2m}^{q-1})$ generates
$\pi_1(L(m,q))$. For we have
\begin{equation*}
F(\xi_{2m}^{q+1},\xi_{2m}^{q-1})(z_0+z_1j)
=\xi_{2m}^{q+1}(z_0+z_1j)\xi_{2m}^{1-q}
=\xi_{m}z_0+\xi_{m}^qz_1j\ .
\end{equation*}

We first collect some information on orientation-reversing isometries.
\begin{proposition} Let $M$ be an elliptic $3$-manifold.
Then $M$ admits an orientation-reversing homeomorphism only if $M$ is a
lens space $L(m,q)$ with $q^2\equiv -1\bmod{m}$. For the lens spaces with
this property:
\begin{enumerate}
\item
The element $T$ of $\Or(4)$ defined by $T(z)=z^{-1}$ induces
orientation-reversing isometries $\tau$ on $S^3$ and on $\R\P^3$, with
$\tau^2=1$.
\item
The element $R\in\Or(4)$ defined by $R(z_0+z_1j)=\overline{z_1}+z_0j$
induces an orientation-reversing isometry $\rho$ on each $L(m,q)$, with
$\rho^2=f(j,j)$ and $\rho^4=1$.
\end{enumerate}
\label{prop:orientation-reversing}
\end{proposition}

\begin{proof}
By~\cite{N-R} (pp.\ 188-190, 194), only the lens spaces $L(m,q)$ with
$q^2\equiv -1\bmod{m}$ can admit orientation-reversing homeomorphisms.  For
$M=S^3$ or $\R\P^3$, $\tau$ clearly has the required properties.  For all
cases with $q^2\equiv -1\bmod{m}$, we have
\begin{gather*}
R\, F(\xi_{2m}^{q+1},\xi_{2m}^{q-1})\,R^{-1}(z_0+z_1j)
=\xi_m^{-q}z_0+\xi_mz_1j\\
=\xi_m^{-q}z_0+\xi_m^{-q^2}z_1j=F(\xi_{2m}^{q+1},\xi_{2m}^q)^{-q}(z_0+z_1j)\
.
\end{gather*}
\noindent Therefore $R$ normalizes $\pi_1(L(m,q))$ and induces an
orientation-reversing isometry $\rho$ on $L(m,q)$. The properties of $\rho$
are immediate from the facts that $R^2=F(j,j)$ and $R^4=1$.
\end{proof}

The next result contains the key group-theoretic facts that determine the
structure of $\Isom(L(m,q))$.

\begin{proposition}
Let $C$ be the cyclic subgroup of order~$m$ in $S^1\ttimes S^1$
generated by $F(\xi_{2m}^{q+1},\xi_{2m}^{q-1})$, where $1<q<m-1$ and
$(m,q)=1$.
\begin{enumerate}
\item
The normalizer of $C$ in $\Ostar\ttimes\Ostar$ always
contains $\Dih(S^1\ttimes S^1)$.
\item
$C$ is normal in $\Ostar\ttimes\Ostar$ if and only if
$q^2\equiv 1\bmod{m}$.
\item
When $C$ is normal in $\Ostar\ttimes\Ostar$, the quotient
$\Ostar\ttimes\Ostar/C$ is isomorphic to the following
groups:
\begin{enumerate}
\item
$\Or(2)\ttimes \Or(2)$ when $(m,q+1)(m,q-1)=m$.
\item
$\Or(2)\times \Or(2)$ when $(m,q+1)(m,q-1)=2m$.
\end{enumerate}
\end{enumerate}
\label{prop:Lie quotients}
\end{proposition}
\noindent In \cite{KM}, these cases are given equivalently as (a) $m$ is
odd or $(q^2-1)/m$ is odd, and (b) both $m$ and $(q^2-1)/m$ are even.

\begin{proof}[Proof of Proposition~\ref{prop:Lie quotients}]
It is clear that $F(j,j)$ normalizes $C$. To determine normality of $C$ in
$\Ostar\ttimes\Ostar$, we need check only one other coset
representative. We have
$$F((j,1)(\xi_{2m}^{q+1},\xi_{2m}^{q-1})(-j,1))(z_0+z_1j)
=\xi_m^{-q}z_0+\xi_m^{-1}z_1j\ ,$$
\noindent so $F(j,1)$ normalizes $C$ if and only if
$$\xi_m^{-q}z_0+\xi_m^{-1}z_1j=F(\xi_{2m}^{q+1},\xi_{2m}^{q-1})^{-q}(z_0+z_1j)
=\xi_m^{-q}z_0+\xi_m^{-q^2}z_1j\ ,$$
\noindent that is, if and only if $q^2\equiv1\bmod{m}$. This establishes
the first two conclusions, so from now on we assume that $C$ is normal.

Let $W$ be the subgroup of $S^1\times S^1$ generated by
$(\xi_{2m}^{q+1},\xi_{2m}^{q-1})$, and let $C^*$ be the preimage of $C$ in
$S^1\times S^1$, that is, the subgroup generated by
$(\xi_{2m}^{q+1},\xi_{2m}^{q-1})$ and $(-1,-1)$.  Of course,
$\Ostar\ttimes\Ostar/C \cong(\Ostar\times\Ostar)/C^*$.  We have
$(\xi_{2m}^{q+1},\xi_{2m}^{q-1})^m=((-1)^{q+1},(-1)^{q-1})$, so when $q$ is
odd, $W$ is cyclic of order~$m$, and $C^*\cong W\times C_2$ with the $C_2$
factor generated by $(-1,-1)$.

Put $d_1=m/(m,q+1)$ and $d_2=m/(m,q-1)$.  When $q$ is odd, $(q+1,q-1)=2$,
and since $(q+1)(q-1)\equiv 0\bmod{m}$, $(m,q+1)(m,q-1)$ is either $m$ or
$2m$, and $d_1d_2$ is either $m$ or~$m/2$. In particular, when $m$ is
odd, $d_1d_2=m$.

Suppose first that $m$ is odd. We may assume that $q$ is odd. For suppose
that $q$ is even. Consider the element $F(\xi_{2m}^{m-q+1},
\xi_{2m}^{m-q-1})$, which generates $\pi_1(L(m,m-q))$. We have
\begin{equation*}
F(\xi_{2m}^{m-q+1},\xi_{2m}^{m-q-1})=
F(-1,-1)F(\xi_{2m}^{-q+1},\xi_{2m}^{-q-1})=
F(\xi_{2m}^{q-1},\xi_{2m}^{q+1})^{-1}\ .
\end{equation*}
\noindent Under the automorphism of $\Ostar\ttimes \Ostar$ which
interchanges the two coordinates, $F(\xi_{2m}^{q-1},\xi_{2m}^{q+1})^{-1}$
is carried to $F(\xi_{2m}^{q+1},\xi_{2m}^{q-1})^{-1}$. Consequently, the
quotients of $\Ostar\ttimes\Ostar$ by the subgroups generated by these two
elements are isomorphic.

Following the notation of~\cite{Orlik}, for a subgroup $H\subseteq
A_1\times A_2$ we will denote by $H_i$ the projections of $H$ to $A_i$ for
$i=1,2$, and put $H_1'=\set{h_1\in H_1\;\vert\;(h_1,1)\in H}$ and
$H_2'=\set{h_2\in H_2\;\vert\;(1,h_2)\in H}$. Then, $H/(H_1'\times
H_2')\cong H_1/H_1'\cong H_2/H_2'$.

Since $q$ is odd, $C^*= W\times C_2$ where $W$ is cyclic of order~$m$
generated by $(\xi_{2m}^{q+1},\xi_{2m}^{q-1})$ and $C_2$ is generated by
$(-1,-1)$.  The order of $W_1$ is $2m/(2m,q+1)=d_1$ and the order of $W_2$
is $2m/(2m,q-1)=d_2$, so the order of $W_1\times W_2$ is $d_1d_2=m$, which
shows that $W=W_1\times W_2\cong C_{d_1}\times C_{d_2}$.  Since $d_1$ and
$d_2$ are odd, the quotient of $\Ostar\times \Ostar$ by $W$ is isomorphic
to $\Ostar\times\Ostar$. The image of $(-1,-1)$ in the quotient is
$(-1,-1)$, so $\Ostar\ttimes\Ostar/C\cong \Ostar\ttimes\Ostar\cong
\Or(2)\ttimes \Or(2)$, as in case~(a).

Now suppose that $m$ is even. Since $q$ is odd, $W$ has order $m$ and
does not contain $(-1,-1)$. Since $W$ must contain an involution, $C^*$
contains the order~$4$ subgroup generated by $(1,-1)$ and $(-1,-1)$. The
quotient of $\Ostar\times \Ostar$ by this order~$4$ subgroup is
$\Or(2)\times \Or(2)$, into which $C^*$ descends to a subgroup $K$ of order
$m/2$ generated by $(\xi_m^{q+1},\xi_m^{q-1})$. The order of $K_1$ is
$m/(m,q+1)=d_1$ and the order of $K_2$ is $m/(m,q-1)=d_2$.

If $d_1d_2=m/2$, then we must have $K=K_1'\times K_2'$, and $\Or(2)\times
\Or(2)/K\cong \Or(2)\times \Or(2)$, giving case~(b). So assume that
$d_1d_2=m$, so that the index of $K_i'$ in $K_i$ is $2$.  The quotient of
$\Or(2)\times \Or(2)$ by $K_1'\times K_2'$ is isomorphic to $\Or(2)\times
\Or(2)$. The image of $K$ must be the subgroup generated by $(-1,-1)$,
since otherwise $K$ would have equaled $K_1'\times K_2'$, and taking the
quotient by this subgroup gives $\Or(2)\ttimes \Or(2)$ as in~(a).
\end{proof}

We now calculate $\Isom(M)$ for the lens spaces.

\medskip
\noindent {\sl Case I.} $m=1$

For $L(1,0)=S^3$, we have $\isom(S^3)=\SO(4)$ and $\Is(S^3)=C_2$ generated
by $\tau$ (from proposition~\ref{prop:orientation-reversing}). Since
$\tau^2=1$ and $\tau\,F(q_1,q_2)\,\tau=F(q_2,q_1)$, $\Isom(S^3)\cong
\Or(4)\cong \SO(4)\circ C_2$ where the $C_2$ acts by interchanging the
factors of $\SO(4)=S^3\ttimes S^3$.
\smallskip

\noindent
In the remaining cases, $m\geq 2$ and we may assume that $1\leq q\leq
m/2$. We have $G_L$ equal to the cyclic group generated by
$\xi_{2m}^{q+1}$, so $G_L$ has order~$2$ if and only $2m=2(q+1)$, that is,
$q=m-1$. Since $q\leq m/2$, this occurs only for $L(2,1)$. On the other
hand, $G_R$ has order $2$ exactly when~$q=1$.

\medskip
\noindent {\sl Case II.} $m=2$

For $L(2,1)=\R\P^3$, we have $G_L\times G_R=C_2\times C_2$, and
$\Norm(G^*)=\SO(4)$, so $\isom(\R\P^3)\cong \SO(3)\times \SO(3)$ and
$\Is(\R\P^3)\cong C_2$ generated by $\tau$. So $\Isom(\R\P^3)\cong
(\SO(3)\times \SO(3))\circ C_2$, where the action of $C_2$ interchanges the
factors.
\smallskip

\noindent {\sl Case III.} $m>2$ and $q=1$.

Since $q^2\not \equiv -1\bmod{m}$, there are no orientation-reversing
isometries.  If $m$ is even, we have $G^*=G_L\times G_R=C_m\times C_2$, so
$\Isom(L(m,1))=(\Ostar\times S^3)/(C_m\times C_2)=\Or(2)\times\SO(3)$. If
$m$ is odd, then $G_L=C_{2m}$ and $G_R=C_2$, and $G^*$ has index $2$ in
$G_L\times G_R$. It is characteristic, being the elements of order less
than $2m$, so $\Norm(G^*)=\Ostar\times S^3$. To find $\Norm(G^*)/G^*$, we
can first take the quotient by the subgroup generated by $(\xi_m,1)$,
obtaining $\Ostar\times S^3$, then take the quotient of this by $(-1,-1)$,
obtaining $\Ostar\ttimes S^3$.  In either case, $\Is(L(m,1))\cong C_2$
generated by $f(j,j)$.
\smallskip

\noindent {\sl Case IV.} $m>2$, $1<q<m/2$ and $q^2\not\equiv\pm1\bmod{m}$
\smallskip

Since $q^2\not \equiv -1\bmod{m}$, there are no orientation-reversing
isometries.  By proposition~\ref{prop:Lie quotients}, $\Norm(G)=
\Dih(S^1\ttimes S^1)$. It follows that $\Isom(L(m,q))\cong \Dih(S^1\ttimes
S^1)\cong \Dih(S^1\times S^1)$, and $\Is(L(m,q))\cong C_2$ generated
by~$f(j,j)$.
\smallskip

\noindent {\sl Case V.} $m>2$, $1<q<m/2$ and $q^2\equiv-1\bmod{m}$
\smallskip

By proposition~\ref{prop:Lie quotients}, $\Norm(G)=\Dih(S^1\ttimes S^1)$,
so $\Is_+(L(m,q))\cong C_2$ generated by $f(j,j)$. The isometry $\rho$ from
proposition~\ref{prop:orientation-reversing} acts on $S^1\ttimes S^1$ by
sending $(z_0,w_0)$ to $(w_0,\overline{z_0})$ and $\rho^2=f(j,j)$, so
$\Isom(L(m,q))\cong (S^1\ttimes S^1)\circ C_4$ where the generator of $C_4$
acts by sending $(z_0,w_0)$ to $(w_0,\overline{z_0})$.
\smallskip

\noindent {\sl Case VI.} $m>2$, $1<q<m/2$ and $q^2\equiv1\bmod{m}$
\smallskip

Since $q^2\not \equiv -1\bmod{m}$, there are no orientation-reversing
isometries.  Pro\-position~\ref{prop:Lie quotients} shows that $\Norm(G)=
\Ostar\ttimes\Ostar$, and gives the quotient groups $\Isom(L(m,q))$. In
both cases, $\Is(L(m,q))\cong C_2\times C_2$ is generated by $f(ij,-i)$,
$f(j,j)$, and $f(i,ij)$.

Note that in all cases other than $M=S^3$ and $M=\R\P^3$,
$\Norm(G)\subseteq \Or(2)\ttimes S^3$.  Table~\ref{tab:lens spaces}
summarizes our calculations of~$\Isom(M)$ for the lens spaces~$L(m,q)$.
The final section of \cite{KM} gives formulas for the number of instances
of each case as a function of $m$.

\begin{table}
\begin{small}
\renewcommand{\arraystretch}{1.5}
\setlength{\tabcolsep}{1.75 ex}
\setlength{\fboxsep}{0pt}
\fbox{%
\begin{tabular}{l|l|l}
$m$, $q$&$\Isom(L(m,q))$&$\Is(L(m,q))$\\
\hline
\hline
$m=1$&$\Or(4)$&$C_2$\\
\hline
$m=2$&$(\SO(3)\times \SO(3))\circ C_2$&$C_2$\\
\hline
\settowidth{\minipagewidth}{$m>2$ even,}%
\begin{minipage}{\minipagewidth}%
\noindent $m>2$ odd,\par\end{minipage} $q=1$&$\Ostar\ttimes S^3$&$C_2$\\
\hline
$m>2$ even, $q=1$&$\Or(2)\times \SO(3)$&$C_2$\\
\hline
$m>2$, $1<q<m/2$, $q^2\not\equiv\pm1\pod{m}$&$\Dih(S^1\times S^1)$&$C_2$\\
\hline
$m>2$, $1<q<m/2$, $q^2\equiv-1\pod{m}$&$(S^1\ttimes S^1)\circ C_4$&$C_4$\\
\hline
$m>2$,
\settowidth{\minipagewidth}{$1<q<m/2$, $q^2\equiv1\pod{m}$,}%
\begin{minipage}[t]{\minipagewidth}%
\noindent $1<q<m/2$, $q^2\equiv1\pod{m}$,\par
\noindent $(m,q+1)(m,q-1)=m$\rule[-1.2 ex]{0mm}{0mm}\par%
\end{minipage}%
&$\Or(2)\ttimes \Or(2)$&$C_2\times C_2$\\
\hline
$m>2$,
\settowidth{\minipagewidth}{$1<q<m/2$, $q^2\equiv1\pod{m}$,}%
\begin{minipage}[t]{\minipagewidth}%
\noindent $1<q<m/2$, $q^2\equiv1\pod{m}$,\par
\noindent $(m,q+1)(m,q-1)=2m$\rule[-1.2 ex]{0mm}{0mm}\par%
\end{minipage}%
&$\Or(2)\times \Or(2)$&$C_2\times C_2$
\end{tabular}}
\end{small}
\bigskip
\caption{Isometry groups of $L(m,q)$}
\label{tab:lens spaces}
\end{table}

\subsection{Realization by finite groups of isometries}
\label{sec:isometric realization}

We have seen that the extension ($*$) is always split.  That is, the
generators of $\Is(M)$ given in sections~\ref{sec:isometries}
and~\ref{sec:lens spaces} actually generate a subgroup of $\Isom(M)$ which
maps isomorphically to $\Is(M)$.  We state this as a ``realization''
theorem.

\begin{theorem} Let $M$ be a an elliptic $3$-manifold.
Then there is a subgroup of $\Isom(M)$ which maps isomorphically to
$\Is(M)$ under the homomorphism $\Isom(M)\to \Is(M)$.
\label{thm:realization}
\end{theorem}

One can check easily that when $\Is(M)$ is nontrivial, the subgroup in
theorem~\ref{thm:realization} is never unique. In fact, there is always a
family of realization subgroups of positive dimension.

\section[Seifert fiberings]
{Seifert fiberings}
\label{sec:fiberings}

In this section, we will examine the extent to which elliptic $3$-manifolds
have Seifert fiberings invariant under their groups of isometries. In
section~\ref{sec:isometries intro}, we imbedded each fundamental group $G$
into $\SO(4)$ in such a way that the Hopf fibering on $S^3$ induces a
Seifert fibering on $S^3/G$, which is maximally symmetric among the
possible Seifert fiberings of $S^3/G$. That is, the subgroup of $\Isom(M)$
that preserves this fibering (i.\ e.\ takes each fiber to a fiber) is as
large as possible. In almost cases, this fibering is preserved by all
isometries.  The precise statement is given in
theorem~\ref{thm:fiber-preserving isometries} below.

We begin by recalling and examining the Hopf fibration of $S^3$, which will
descend to the maximally symmetric fibration on $S^3/G$. For the rest of
this section, we regard the $2$-sphere $S^2$ as $\C\cup\set{\infty}$.  We
speak of antipodal points and orthogonal transformations on $S^2$ by
transferring them from the standard $2$-sphere using the stereographic
projection. In particular, the antipodal map is defined by
$\alpha(z)=-1/\overline{z}$.

As is well-known, the Hopf fibering on $S^3$ is an $S^1$-bundle structure
with projection map $H\colon S^3\to S^2$ defined by
$H(z_0,z_1)=z_0/z_1$. The fibers are the orbits of the left action of $S^1$
on $S^3$.  The element $F(j,1)$ preserves these fibers and induces the
antipodal map on $S^2$.  For we have
$j(z_0+z_1j)=-\overline{z_1}+\overline{z_0}$, so
$H(F(j,1)(z_0+z_1j))=-1\,/\,\overline{z_0/z_1}$.  Since right
multiplication by elements of $S^3$ commutes with the left action of $S^1$,
it preserves the Hopf fibering, and there is an induced action of $S^3$ on
$S^2$. In fact, it acts orthogonally. For if we write $x=x_0+x_1j$ and
$z=z_0+z_1j$, we have
$F(1,x)(z)=zx^{-1}=z_0\overline{x_0}+z_1\overline{x_1} +(z_1x_0-z_0x_1)j$,
so the induced action on $S^2$ carries $z_0/z_1$ to $
(z_0\overline{x_0}+z_1\overline{x_1})/(z_1x_0-z_0x_1)
=\begin{pmatrix}\phantom{-}\overline{x_0}&\overline{x_1}\\
-x_1&x_0\end{pmatrix}(z_0/z_1)$.  The trace of this linear fractional
transformation is real and lies between $-2$ and $2$ (unless $x=\pm 1$,
which acts as the identity on $S^2$), so it is elliptic.  Its fixed points
are $\big(\,(x_0-\overline{x_0})\pm\sqrt{(x_0-\overline{x_0})^2 -
4x_1\overline{x_1}\,}\,\,\big)/(2x_1)$, which are antipodal, so it is an
orthogonal transformation. Combining these observations, we see that the
action induced on $S^2$ via $H$ determines a surjective homomorphism
$h\colon \Ostar\times S^3\to\Or(3)$, given by $h(x_0,1)=1$ for $x_0\in
S^1$, $h(j,1)$ is the antipodal map, and $h(1,x_0+x_1j)=
\begin{pmatrix}\phantom{-}\overline{x_0}&\overline{x_1}\\
-x_1&x_0\end{pmatrix}$. The kernel of $h$ is $S^1\ttimes\set{\pm 1}\cong
S^1$.

Since the $G$ constructed above all lie in $\Ostar\ttimes S^3$, they
preserve the Hopf fibering on $S^3$, which descends to a fibering on $M$
that we call the \emph{Hopf fibering} on $M$. Its quotient orbifold $\orb$
is $S^2/h(G)$, the quotient of $S^2$ by the action of $h(G)$, and can
easily be calculated using our descriptions of $\Isom(M)$. Since we have
selected our fiberings to be maximally symmetric, the quotient orbifolds in
the cases of $G=Q_8$ and $G=D_{4m}$ are perhaps not the most ``standard''
ones. Table~\ref{tab:orbifolds} lists the quotient orbifolds $\orb$ for all
cases. There, $S^2(\alpha_1,\dots,\alpha_r)$ denotes the orbifold with
underlying topological space the $2$-sphere $S^2$ and cone points of orders
$\alpha_1,\dots\,$, $\alpha_r$, and $(\R\P^2;)$ denotes $\R\P^2$ regarded
as an orbifold.  These orbifolds all carry metrics of constant positive
curvature inherited from $S^2$. Table~\ref{tab:orbifolds} also gives
$d(\orb)$, the ``degree of symmetry'', which is the maximal dimension for a
compact Lie group acting effectively as orbifold diffeomorphisms on~$\orb$.

\begin{table}
\begin{small}
\renewcommand{\arraystretch}{1.5}
\setlength{\tabcolsep}{1.5 ex}
\setlength{\fboxsep}{0pt}
\fbox{%
\begin{tabular}{l|l|l|c}
$G$&$M$&$\orb=S^2/h(G)$&$d(\orb)$\\ 
\hline 
\hline
$\set{1}$&$S^3$&$(S^2;)$&$3$\\ 
\hline
$C_2$&$\R\P^3$&$(\R\P^2;)$&$3$\\
\hline
$C_m$&$L(m,1)$&$(S^2;)$&$3$\\
\hline 
$C_m$&$L(m,q)$, $1<q<m/2$&%
\settowidth{\minipagewidth}{$\;\;k=m/(q-1,m)$}%
\begin{minipage}{\minipagewidth}{$(S^2;k,k),$\rule[2.8 ex]{0mm}{0mm}\par%
$\;\;k=m/(q-1,m)$\rule[-1.2 ex]{0mm}{0mm}\par}%
\end{minipage}%
&$1$\\ 
\hline 
$Q_8$&quaternionic manifold&$(\R\P^2;)$&$3$\\ 
\hline 
$Q_8\times C_n$&quaternionic manifold&$(S^2;2,2,2)$&$0$\\ 
\hline 
$D_{4m}^*$&prism manifold&$(\R\P^2;)$&$3$\\ 
\hline 
$D_{4m}^*\times C_n$&prism manifold&$(S^2;2,2,m)$&$0$\\ 
\hline 
index $2$ diagonal&prism manifold&$(S^2;2,2,m)$&$0$\\ 
\hline 
&tetrahedral manifold&$(S^2;2,3,3)$&$0$\\ 
\hline 
&octahedral manifold&$(S^2;2,3,4)$&$0$\\ 
\hline 
&icosahedral manifold&$(S^2;2,3,5)$&$0$\\
\end{tabular}}
\end{small}
\bigskip
\caption{Quotient orbifolds for the Hopf fiberings}
\label{tab:orbifolds}
\end{table}

The complete sets of isomorphism classes of Seifert fiberings for elliptic
$3$-manifolds are detailed in~\cite{Orlik} and~\cite{Sakuma}, the latter
including corrections of a couple of minor errors in the former. The lens
spaces have infinitely many fiberings, but for the other cases, the
possible fiberings can be understood quite simply. Suppose that $M$ is a
quaternionic or prism manifold. In this case, there are two isomorphism
classes of fiberings corresponding to the fact that $G$ (as we have
imbedded it in $\SO(4)$) preserves both the Hopf fibering on $S^3$ and the
mirror Hopf fibering, which is the image of the standard Hopf fibering on
$S^3$ under the involution of taking the quaternionic inverse (its fibers
are the orbits of the action of $S^1$ by right multiplication).  Depending
on the imbedding of $G$, one of the Hopf or mirror Hopf fiberings on $M$
has a quotient orbifold with underlying topological space $\R\P^2$, and the
other has quotient orbifold $(S^2;2,2,2)$ or $(S^2;2,2,m)$. Finally, when
$M$ is a tetrahedral, octahedral, or icosahedral manifold, the mirror Hopf
fibering of $S^3$ is not invariant, and the unique fibering on $M$ is the
Hopf fibering. The extent to which $\Isom(M)$ preserves the fiberings is
detailed in the next result.

\begin{theorem}
Let $G$ be one of the subgroups of $\SO(4)$ acting freely on $S^3$ and let
$M=S^3/G$ with the Hopf fibering.
\begin{enumerate}
\item[(a)] The connected component of the identity $\isom(M)$ preserves the
Hopf fibering in all cases except when $G$ is one of $\set{1}$, $C_2$,
$T_{24}^*$, $O_{48}^*$, or $I_{120}^*$. For these exceptional manifolds, no
Seifert fibering is preserved by~$\isom(M)$.
\item[(b)] The full orientation-preserving isometry group $\Isom_+(M)$
preserves the Hopf fibering in all cases except those in \textup{(a)},
together with $G=Q_8$. For these exceptional manifolds, no Seifert fibering
is preserved by~$\Isom_+(M)$.
\item[(c)] No orientation-reversing isometry of any $M$ can preserve a
Seifert fibering.
\end{enumerate}
\label{thm:fiber-preserving isometries}
\end{theorem}

\begin{proof}
The fact that the isometries preserve the fiberings in the nonexceptional
cases in~(a) and~(b) is immediate from fact that $\norm(G)$ or $\Norm(G)$
is contained in $\Ostar\ttimes S^3$, as remarked at the ends of
section~\ref{sec:isometries} and~\ref{sec:lens spaces}. 

For the exceptional cases in~(a), the dimension of the group of isometries
is~$6$ for the cases of $S^3$ and $\R\P^3$ and is $3$ for the other
cases. Since the isometries act as a group of orbifold diffeomorphisms on
the quotient orbifold, the dimension of a group of fiber-preserving
isometries is at most $1$ more than the degree of symmetry of the quotient
orbifold, the extra dimension coming from the isometries that move fibers
vertically. For $S^3$ and $\R\P^3$ this sum is $4$, and for the other cases
in~(a) it is $1$, so $\isom(M)$ cannot preserve the Hopf fibering.  For
$S^3$ and $\R\P^3$, there are other fiberings, but none can be preserved by
a group of dimension greater than~$4$, since the maximum degree of symmetry
of a quotient orbifold is $3$ (achieved by the Hopf fiberings).  For the
other three cases in (a), the fiberings are unique up to isomorphism and
hence all fiberings have the same quotient orbifold, whose degree of
symmetry is $0$. Since $\isom(M)$ has dimension $3$, it cannot preserve any
fibering.

For the additional quaternionic manifold in~(b), some of the isometries
$f(q,1)$ with $q\in O_{48}^*-Q_8$ do not preserve the Hopf fibering. In
fact, $\Isom_+(M)$ cannot preserve any Seifert fibering isomorphic to this
one. For the fibering is nonsingular, with homotopy exact sequence $1\to
C_4\to \pi_1(M)\to C_2\to 1$, where $C_4$ is the subgroup of $\pi_1(M)$
represented by the fiber. Since no cyclic subgroup of order~$4$ of $Q_8$ is
invariant under all automorphisms induced by elements of $\Isom_+(M)$, no
such fibering can be invariant. The quotient orbifold for the other Seifert
fibering on this manifold is $(S^2;2,2,2)$. Some isometries from each path
component do preserve this fibering, but since $d(S^2;2,2,2)=0$, these
other fiberings cannot be preserved by all of~$\Isom_+(M)$.

Observation~(c) is due to the fact that the Seifert fiberings on the $M$
all have nonzero Euler class. Any fiber-preserving diffeomorphism must
preserve the Euler class, but any orientation-reversing diffeomorphism must
multiply it by~$-1$, so there are no fiber-preserving orientation-reversing
isometries.
\end{proof}

\section[Mapping class groups]
{Mapping class groups}
\label{sec:mapping class groups}

Through the efforts of a number of authors, the mapping class groups of all
elliptic $3$-manifolds have been determined. Using these results, we will
verify the ``$\pi_0$'' part the Generalized Smale Conjecture, as discussed
in the introduction. The statement is as follows, where $\H(M)$ to denote
the (not necessarily orientation-preserving) mapping class
group~$\pi_0(\Diff(M))$.

\begin{theorem}
Let $M$ be an elliptic $3$-manifold. Then the inclusion of $\Isom(M)$ into
$\Diff(M)$ is a bijection on path components. That is, it induces an
isomorphism from $\Is(M)$ to~$\H(M)$.
\label{thm:pi0 Smale}
\end{theorem}

\begin{proof}
For $S^3$ this was first proven by Cerf (see~\cite{C}). For the remaining
cases, we will rely on calculations of $\H(M)$ due to several authors.  We
utilize the composite homomorphism $\Phi\colon \Is(M)\to \H(M)\to\Out(G)$
which takes each isometry to the outer automorphism it induces on
$\pi_1(M)=G$. As noted in section~\ref{sec:calculation of isometries},
elementary covering space theory shows that $\Phi(f(q_1,q_2))$ is the outer
automorphism induced on $G$ by conjugation by $F(q_1,q_2)$.

In all cases except the tetrahedral and icosahedral manifolds, there are
explicit published calculations of $\H(M)$ which show it to be abstractly
isomorphic to $\Is(M)$. In these cases, it is sufficient to check that
$\Phi$ is injective.

For lens spaces other than $S^3$, Bonahon~\cite{B} gave a complete
calculation of $\H(L(m,q))$, obtaining groups isomorphic to the
$\Is(L(m,q))$ obtained in section~\ref{sec:lens spaces} (these were also
calculated in \cite{HR}). In particular, $\H(\R\P^3)\cong C_2$, and since
$\Is(\R\P^3)$ contains orientation-reversing elements,
$\Is(\R\P^3)\to\H(\R\P^3)$ is an isomorphism. In all cases with $m>2$,
$\Phi$ is injective. For if we regard our standard generator
$F(\xi_{2m}^{q+1},\xi_{2m}^{q-1})$ of $\pi_1(L(m,q))$ as $t\in \langle
t\;\vert\;t^m=1\rangle$, we can check that $f(j,j)$ sends $t$ to $t^{-1}$,
$f(j,1)$ sends $t$ to $t^q$, and $\rho$ sends $t$ to $t^q$, an automorphism
of order~$4$ in the cases when $q^2\equiv -1\bmod{m}$. In all cases $\Phi$
is injective.

For the $M$ that contain one-sided Klein bottles, the mapping class groups
were calculated in~\cite{A} and~\cite{R2}. These include all quaternionic
and prism manifolds, and the lens spaces $L(4k,2k-1)$. Again, the groups
are isomorphic to~$\Is(M)$.

For the quaternionic manifolds, the abelianization $\Ab(Q_8\times C_n)$ is
isomorphic to $(C_2\times C_2)\times C_n$. The $S_3$ factors in the
isometry groups act as $\GL(2,2)$ on the $C_2\times C_2$ factor. When
$n\neq 1$, the isometry $f(j,1)$ acts as multiplication by $-1$ on the
$C_n$ factor.  Since $n$ is odd, we have $n>2$, so this is a nontrivial
automorphism on $\Ab(G)$.  It follows that $\Phi$ is injective for these
cases.

For the prism manifolds, note first that $f(1,\xi_{4m})$ sends $j$ to
$\xi_{2m}j$, which is the nontrivial outer automorphism of $D_{4m}^*$.  For
the cases of $D_{4m}^*\times C_n$, $n$ is odd, and $f(j,1)$ induces
multiplication by $-1$ on the $C_n$ factor, so $\Phi$ is injective for
these cases as well. When $G$ is an index~$2$ diagonal subgroup in
$C_{4n}\times D_{4m}^*$, the automorphism induced by the nontrivial element
of $\Is(M)$ multiplies by $-1$ in the $C_{4n}$ factor, so is not inner.

The mapping class groups of the octahedral manifolds were calculated
in~\cite{R-B}, obtaining groups isomorphic to $\Is(M)$. The only nontrivial
element is $f(j,1)$, for the cases $O_{48}^*\times C_n$. Since $n>2$, this
induces a nontrivial outer automorphism.

For the remaining cases, we do not know an explicit calculation of
$\H(M)$. But as the authors of~\cite{B-O} note (and carry out for the case
of $G=I_{120}^*$, for which $S^3/G$ is the Poincar\'e sphere), the
following result from their paper leads easily to one.

\begin{theorem}[Boileau-Otal]
Let $M$ be a tetrahedral or icosahedral manifold. Then any diffeomorphism
of $M$ is isotopic to a fiber-preserving diffeomorphism. Consequently, any
diffeomorphism which is homotopic to the identity is isotopic to the
identity map.
\label{thm:Boileau-Otal}
\end{theorem}

\noindent To utilize this theorem, we use the following lemma.
\begin{lemma}
Let $M$ be an elliptic $3$-manifold $S^3/G$ with $\abs{G}>2$. Suppose that
$f\colon M\to M$ is a diffeomorphism which induces the identity outer
automorphism on $\pi_1(M)=G$. Then $f$ is homotopic to the identity map.
\label{lem:homotopy}
\end{lemma}

\begin{proof}
Triangulate $M$. We may assume that $f$ is the identity on the $1$-skeleton
of $M$ (first make it preserve the basepoint and induce the identity
automorphism, then make it fix a maximal tree in the $1$-skeleton
$M^{(1)}$, then since $f$ induces the identity automorphism we may deform
it to be the identity on the remaining edges of $M^{(1)}$). Since
$\pi_2(M)=\set{0}$, we may further make $f$ the identity on
$M^{(2)}$. There is a subcomplex $K$ of $M^{(2)}$ such that $M$ is obtained
by attaching a single $3$-cell $C$ to $K$ along $\partial C$. The
characteristic map of $C$ and its composition with $f$ agree on $\partial
C$, so taking them as the maps on the top and bottom hemispheres defines a
map $g\colon S^3\to M$. Via the isomorphism $H^3(M,K;\Z)\cong H^3(M;\Z)$,
the degree $\deg(f)$ equals $1+\deg(g)$. Since $g$ factors through the
universal covering $S^3\to M$, $\deg(g)=k\abs{G}$ for some $k$. Since
$\deg(f)=\pm1$ and $\abs{G}>2$, it follows that $g$ has degree $0$ and
hence $f$ is homotopic to the identity map.
\end{proof}

\noindent By the lemma, theorem~\ref{thm:Boileau-Otal} shows that
$\H(M)\to\Out(G)$ is injective.

For the remaining cases, note that if $n$ is relatively prime to the order
of the finite group $A$, then $\Out(A\times C_n)\cong \Out(A)\times
\Aut(C_n)$. In each of these cases, we will show that $\Phi$ is injective
and that all outer automorphisms of $G$ that are realizable by
diffeomorphisms lie in the image of $\Phi$, from which it follows using the
injectivity of $\H(M)\to\Out(G)$ that $\Is(M)\to\H(M)$ is an isomorphism.

Suppose first that $G=T_{24}^*$. The abelianization $\Ab(T_{24})^*$ is
$C_3$, on which $f(1,\frac{1}{\sqrt{2}}i+\frac{1}{\sqrt{2}}j)$ induces the
nontrivial automorphism.  From~\cite{Plotnick}, $\Out(T_{24}^*)$ is $C_2$,
so $\Phi$ is an isomorphism.

Suppose that $G=T_{24}^*\times C_n$ with $n>1$ (and hence $n>2$ since $n$
must be odd).  By theorem~\ref{thm:Boileau-Otal}, every diffeomorphism is
isotopic to one that preserves the fiber, so must induce multiplication by
$1$ or $-1$ on the $C_n$ factor which is represented by the fiber. So the
image of $\H(M)\to \Out(G)$ must be contained in $C_2\times
\Out(T_{24}^*)\cong C_2\times C_2$, and $\Phi$ carries $\Is(M)$
isomorphically to this subgroup.

Finally, suppose that $G$ is an index~$3$ diagonal subgroup in
$C_{6n}\times T_{24}^*$. As before, any diffeomorphism must induce
multiplication by $\pm1$ in the first factor. Multiplication by $-1$ only
preserves the diagonal subgroup if we also apply the nontrivial outer
automorphism to the elements of the second factor. So the automorphisms
realizable by diffeomorphisms lie in a $C_2$ subgroup of $\Out(G)$, and
$\Phi$ carries $\Is(M)$ isomorphically to this subgroup.

For $G\cong I_{120}^*$, we have from~\cite{Plotnick} that
$\Out(I_{120}^*)\cong C_2$, but that the nontrivial element is not induced
by any diffeomorphism or even homotopy equivalence of $M$. So (as
calculated explicitly in~\cite{B-O}), $\H(M)=\set{1}$, and we have seen
that $\Is(M)=\set{1}$ as well. For $G\cong I_{120}^*\times C_n$, the
nontrivial element of $\Out(I_{120}^*)$ cannot be in the image of
$\H(M)\to\Out(G)$, for if it could be realized by a diffeomorphism, then
this would lift to a diffeomorphism of the covering corresponding to the
subgroup $I_{120}^*\times \set{1}$ (this subgroup is characteristic since
$(30,n)=1$). The lift would induce the nontrivial outer automorphism, and
we have already seen that this is impossible. So the only realizable outer
automorphism is multiplication by $-1$ on $C_n$, and $\Phi$ carries
$\Is(M)$ isomorphically to this subgroup.
\end{proof}

Combining theorem~\ref{thm:pi0 Smale} with theorem~\ref{thm:realization},
we have a Nielsen realization theorem for elliptic $3$-manifolds.

\begin{theorem}
Let $M$ be an elliptic $3$-manifold. Then there is a subgroup of $\Isom(M)$
that maps isomorphically to $\H(M)$ under the natural homomorphism
$\Isom(M)\to \H(M)$.
\label{thm:Nielsen realization}
\end{theorem}

\bibliographystyle{amsplain}

\end{document}